\documentclass[11pt]{amsart}
\usepackage{enumerate}
\usepackage{amsfonts}
\usepackage{amssymb}

\makeatletter
\@namedef{subjclassname@2010}{%
  \textup{2010} Mathematics Subject Classification}
\makeatother

\newtheorem{theorem}{Theorem}

\newtheorem{proposition}{Proposition}
\newtheorem{lemma}{Lemma}
\newtheorem{claim}{Claim}

\theoremstyle{definition}
\newtheorem{definition}{Definition}
\newtheorem{problem}{Problem}

\input xy
\xyoption{all}

\newcommand{\U}{\mathcal U}
\newcommand{\V}{\mathcal V}
\newcommand{\W}{\mathcal W}

\newcommand{\ord}{\mathrm{ord}}
\newcommand{\F}{\mathcal{F}}

\newcommand{\w}{\omega}
\newcommand{\IN}{\mathbb N}
\newcommand{\IC}{\mathbb C}
\newcommand{\hT}{\widehat{\mathsf K}}

\newcommand{\pr}{\operatorname{pr}}
\newcommand{\St}{\mathcal S{t}}

\begin{document}

\title{Characterizing chainable, tree-like, and circle-like continua}
\author{Taras Banakh, Zdzis{\l}aw Koszto{\l}owicz, S{\l}awomir Turek}
\address[T.Banakh, Z.Koszto\l owicz, S.Turek]{Institute of Mathematics, Uniwersytet Humanistyczno-Przyrodniczy
Jana Kochanowskiego, ul. \'{S}wietokrzyska~15, 25-406 Kielce, Poland}
\email{zdzisko@ujk.kielce.pl, sturek@ujk.kielce.pl}
\address[T.Banakh]{Department of Mathematics, Ivan Franko National University of Lviv, Ukraine}
\email{tbanakh@yahoo.com}

\date{}

\begin{abstract} We prove that a continuum $X$ is tree-like (resp. circle-like, chainable) if and only if for each open cover $\U_4=\{U_1,U_2,U_3,U_4\}$ of $X$ there is a $\U_4$-map $f:X\to Y$ onto a tree (resp. onto the circle, onto the interval). A continuum $X$ is an acyclic curve if and only if for each open cover $\U_3=\{U_1,U_2,U_3\}$ of $X$ there is a $\U_3$-map $f:X\to Y$ onto a tree (or the interval $[0,1]$).
\end{abstract}



\maketitle

\section{Main results}
In this paper we characterize chainable, tree-like and circle-like continua in the spirit of the following Hemmingsen's  characterization of covering dimension \cite[1.6.9]{En}.

\begin{theorem}[Hemmingsen]\label{hemming} For a compact Hausdorff space $X$ the following conditions are equivalent:
\begin{enumerate}
\item $\dim X\le n$, which means that any open cover $\U$ of $X$ has an open refinement $\V$ of order $\le n+1$;
\item each open cover $\U$ of $X$ with cardinality $|\U|\le n+2$ has an open refinement $\V$ of order $\le n+1$;
\item each open cover $\{U_i\}_{i=1}^{n+2}$ of $X$ has an open refinement $\{V_i\}_{i=1}^{n+2}$ with $\bigcap_{i=1}^{n+2}V_i=\emptyset$.
\end{enumerate}
\end{theorem}

We say that a cover $\V$ of $\U$ is a {\em refinement} of a cover $\U$ if each set $V\in\V$ lies in some set $U\in\U$. The {\em order} of a cover $\U$ is defined as the cardinal 
$$\ord(\U)=\sup\{|\F|:\F\subseteq\U\mbox{ with $\bigcap$}\F\ne\emptyset\}.$$

An open cover $\U$ of $X$ is called
\begin{itemize}
\item a {\em chain-like} if for $\U$ there is an enumeration $\U=\{U_1,\dots,U_n\}$ such that $U_i\cap U_j\ne\emptyset$ if and only if $|i-j|\le1$ for all $1\le i,j\le n$;
\item {\em circle-like} if there is an enumeration $\U=\{U_1,\dots,U_n\}$ such that $U_i\cap U_j\ne\emptyset$ if and only if $|i-j|\le 1$ or $\{i,j\}=\{1,n\}$; 
\item a {\em tree-like} if $\U$ contains no circle-like subfamily $\V\subseteq\U$ of cardinality $|\V|\ge 3$.
\end{itemize}

We recall that a continuum $X$ is called {\em chainable} (resp. {\em tree-like}, {\em circle-like}) if each open cover of $X$ has a chain-like (resp. tree-like, circle-like) open refinement. By a {\em continuum} we understand a connected compact Hausdorff space.

The following characterization of chainable, tree-like and circle-like continua is the main result of this paper. For chainable and tree-like continua this characterization was announced (but not proved) in \cite{BB}.

\begin{theorem}\label{main} A continuum $X$ is chainable (resp. tree-like, circle-like) if and only if any open cover $\U$ of $X$ of cardinality $|\U|\le 4$ has a chain-like (resp. tree-like, circle-like) open refinement.
\end{theorem} 

In fact, this theorem will be derived from a more general theorem treating $\mathsf K$-like continua.

\begin{definition} Let $\mathsf K$ be a class of continua and $n$ be a cardinal number. A continuum $X$ is called {\em $\mathsf K$-like} (resp. $n$-$\mathsf K$-{\em like}) if for any open cover $\U$ of $X$ (of cardinality $|\U|\le n$) there is a $\U$-map $f\colon X\to K$ onto some space $K\in\mathsf K$.
\end{definition}

We recall that a map $f:X\to Y$ between two topological spaces is called a {\em $\U$-map}, where $\U$ is an open cover of $X$, if there is an open cover $\V$ of $Y$ such that the cover $f^{-1}(\V)=\{f^{-1}(V):V\in\V\}$ refines the cover $\U$. It worth mentioning that a closed map $f:X\to Y$ is a $\U$-map if and only if the family $\{f^{-1}(y)\colon y\in Y\}$ refines $\U$.  

It is clear that a continuum $X$ is tree-like  (resp. chainable, circle-like) if and only if it is $\mathsf K$-like for the class $\mathsf K$ of all trees (resp. for $\mathsf K=\{[0,1]\}$, $\mathsf K=\{S^1\}$). Here $S^1=\{z\in\IC:|z|=1\}$ stands for the circle. 

It turns out that each 4-$\mathsf K$-like continuum is $\hT$-like for some extension $\hT$ of the class $\mathsf K$. This extension is defined with help of locally injective maps.

A map $f:X\to Y$ between topological spaces is called {\em locally injective} if each point $x\in X$ has a neighborhood $O(x)\subseteq X$ such that the restriction $f\restriction O(x)$ is injective. For a class of continua $\mathsf K$ let $\hT$ be the class of all continua $X$ that admit a locally injective map $f:X\to Y$ onto some continuum $Y\in\mathsf K$.

\begin{theorem}\label{hat} Let $\mathsf K$ be a class of 1-dimensional continua. If a continuum $X$ is $4$-$\mathsf K$-like, then $X$ is $\hT$-like.
\end{theorem}

 In Proposition~\ref{li} we shall prove that each locally injective map $f:X\to Y$ from a continuum $X$ onto a tree-like continuum $Y$ is a homeomorphism. This implies that $\hT=\mathsf K$ for any class $\mathsf K$ of tree-like continua. 
This fact combined with Theorem~\ref{hat} implies the following characterization:

\begin{theorem}\label{Th3} 
Let $\mathsf K$ be a class of tree-like continua. A continuum $X$ is $\mathsf K$-like if and only if it is $4$-$\mathsf K$-like.
\end{theorem}

One may ask if the number 4 in this theorem can be lowered to 3 as in the Hemmingsen's characterization of 1-dimensional compacta. It turns out that this cannot be done: the $3$-$\mathsf K$-likeness is equivalent to being an acyclic curve. A continuum $X$ is called a {\em curve} if $\dim X\le 1$. It is {\em acyclic} if each map $f:X\to S^1$ to the circle is null-homotopic.

\begin{theorem}\label{t3} Let $\mathsf K\ni[0,1]$ be a class of tree-like continua. A continuum $X$ is $3$-$\mathsf K$-like if and only if $X$ is an acyclic curve.
\end{theorem}

It is known that each tree-like continuum is an acyclic curve but there are acyclic curves, which are not tree-like \cite{CC}. On the other hand, each locally connected acyclic curve is tree-like (moreover, it is a dendrite \cite[Chapter X]{Nad}).
Therefore, for any continuum $X$ and a class $\mathsf K\ni[0,1]$ of tree-like continua we get the following chain of equivalences and implications (in which the dotted implication holds under the additional assumption that the continuum $X$ is locally connected):
$$
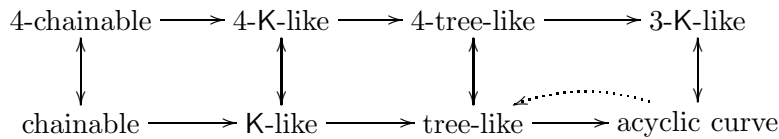
\xymatrix{
\mbox{$4$-chainable}\ar[r]\ar@{<->}[d]&\mbox{$4$-$\mathsf K$-like}\ar@{<->}[d]\ar[r]&\mbox{$4$-tree-like}\ar@{<->}[d]\ar[r]&\mbox{$3$-$\mathsf K$-like}\ar@{<->}[d]\\
\mbox{chainable}\ar[r]&\mbox{$\mathsf K$-like}\ar[r]&\mbox{tree-like}\ar[r]&\mbox{acyclic curve}\ar@{.>}@/_1pc/[l]
}$$

Finally, let us present a factorization theorem that reduces the problem of studying $n$-$\mathsf K$-like continua to the metrizable case.
It will play an important role in the proof of the ``circle-like'' part of Theorem~\ref{main}. 

\begin{theorem}\label{metr} Let $n\in\IN\cup\{\w\}$ and $\mathsf K$ be a family of metrizable continua. A continuum $X$ is $n$-$\mathsf K$-like if and only if any map $f:X\to Y$ to a metrizable compact space $Y$ can be written as the composition $f=g\circ \pi$ of a continuous map $\pi:X\to Z$ onto a metrizable $n$-$\mathsf K$-like continuum $Z$ and a continuous map $g:Z\to Y$.
\end{theorem}

\section{Proof of Theorem~\ref{t3}} Let $\mathsf K\ni[0,1]$ be a class of tree-like continua. We need to prove that a continuum $X$ is $3$-$\mathsf K$-like if and only if $X$ is an acyclic curve.

To prove the ``if'' part, assume that $X$ is an acyclic curve. By Theorem~2.1 of \cite{BB}, $X$ is 3-chainable. Since $[0,1]\in\mathsf K$, the continuum $X$ is $3$-$\mathsf K$-like and we are done.

Now assume conversely, that a continuum $X$ is $3$-$\mathsf K$-like. First, using  Hemmingsen's Theorem~\ref{hemming}, we shall show that $\dim X\le 1$. 
Let $\V=\{V_1,V_2,V_3\}$ be an open cover of $X$. Since the space $X$ is $3$-$\mathsf{K}$-like, we can find a $\V$-map $f:X\to T$ onto a tree-like continuum $T$. Using the 1-dimensionality of tree-like continua, find an open cover $\W$ of $T$ order $\le 2$ such that the cover $f^{-1}(\W)=\{f^{-1}(W):W\in\W\}$ is a refinement of $\V$. 
The continuum $X$ is $1$-dimensional by the implication (2)$\Rightarrow$(1) of Hemmingsen's theorem.

It remains to prove that $X$ is acyclic. Let $f\colon X\to S^1$ be a continuous map. Let $\mathcal{U}=\{U_1,U_2,U_3\}$ be a 
cover of the unit circle $S^1=\{z\in\IC:|z|=1\}$ by three open arcs $U_1,U_2,U_3$, each of length $<\pi$. Such a cover necessarily has $\ord(\mathcal{U})=2$.
By our assumption there is an open finite cover $\mathcal{V}$ of $X$ inscribed in $\{f^{-1}(U_i)\colon i=1,2,3\}$. So, there is a tree-like continuum $T\in\mathsf{K}$ and $\mathcal{V}$-map $g\colon X\to T$.
We can assume that $T$ is a tree and $\mathcal{V}$ is a tree-open cover of $X$. It is well known (see e.g.~\cite{CC}) that there exists a continuous map $h\colon T\to S^1$ that $h\circ g$ is homotopic to $f$. But each map from a tree to the circle is null-homotopic. Hence $h\circ g$ as well $f$ is null-homotopic too.

\section{Proof of Theorem~\ref{hat}}

We shall use some terminology from Graph Theory. So at first we recall some definitions.

By a ({\em combinatorial\/}) {\em graph} we understand a pair $G=(V,E)$ consisting of a finite set $V$ of vertices and a set $E\subseteq\big\{\{a,b\}:a,b\in V,\; a\ne b\big\}$ of unordered pairs of vertices, called {\em edges}. A graph $G=(V,E)$ is {\em connected} if any two distinct vertices $u,v\in V$ can be linked by a path $(v_0,v_1,\dots,v_n)$ with $v_0=u$, $v_n=v$ and $\{v_{i-1},v_i\}\in E$ for $i\le n$. The number $n$ is called the length of the path (and is equal to the number of edges involved). Each connected graph possesses a natural path metric on the set of vertices $V$: the distance between two distinct vertices equals the smallest length of a path linking these two vertices.

Two vertices $u,v\in V$ of a graph are {\em adjacent} if $\{u,v\}\in E$ is an edge. The {\em degree} $\deg(v)$ of a vertex $v\in V$ is the number of vertices $u\in V$ adjacent to $v$ in the graph. The number $\deg(G)=\max\limits_{v\in V}\deg(v)$ is called the {\em degree} of a graph. By an {\em $r$-coloring of a graph} we understand any map $\chi:V\to r=\{0,\dots,r-1\}$. In this case for a vertex $v\in V$ the value $\chi(v)$ is called the {\em color} of $v$.

\begin{lemma}\label{coloring} Let $G=(V,E)$ be a connected graph with $\deg(G)\le 3$ such that $d(u,v)\ge6$ for any two vertices $u,v\in V$ of order 3. Then there is a 4-coloring $\chi\colon V\to4$ such that no distinct vertices $u,v\in V$ with $d(u,v)\le2$ have the same color.
\end{lemma}

\begin{proof} Let $V_3=\{v\in V:\deg(v)=3\}$ denote the set of vertices of order 3 in $G$ and let $\bar B(v)=\{v\}\cup\{u\in V\colon \{u,v\}\in E\}$ be the unit ball centered at $v\in V$.  It follows from $\deg(G)\le 3$ that $|\bar B(v)|\le 4$ for each $v\in V$. Moreover, for any distinct vertices $v,u\in V_3$ the balls $\bar B(v)$ and $\bar B(u)$ are disjoint (because $d(v,u)\ge 6>2$). 
So we can define a 4-coloring $\chi$ on the union $\bigcup_{v\in V_3}\bar B(v)$ so that $\chi$ is injective on each ball $\bar B(u)$ and $\chi(v)= \chi(w)$ for each $v,w\in V_3$. Next, it remains to color the remaining vertices all of order $\le 2$ by four colors so that $\chi(x)\ne\chi(y)$ if $d(x,y)\le 2$. It is easy to check that this always can be done. 
\end{proof}

Each graph $G=(V,E)$ can be also thought as a topological object: just embed the set of vertices $V$ as a linearly independent subset into a suitable Euclidean space  and consider the union $|G|=\bigcup_{\{u,v\}\in E}[u,v]$ of intervals corresponding to the edges of $G$. Assuming that each interval $[u,v]\subseteq |G|$ is isometric to the unit interval $[0,1]$, we can extend the path-metric of $G$ to the path-metric $d$ on the geometric realization $|G|$ of $G$. For a point $x\in|G|$ by $B(x)=\{y\in|G|:d(x,y)<1\}$ and $\bar B(x)=\{y\in|G|:d(x,y)\le 1\}$ denote respectively the open and closed unit balls centered at $x$. More generally,  by $B_r(x)=\{y\in|G|:d(x,y)<r\}$ we shall denote the open ball of radius $r$ with center at $x$ in $|G|$.

By a {\em topological graph} we shall understand a topological space $\Gamma$ homeomorphic to the geometric realization $|G|$ of some combinatorial graph $G$. 
In this case $G$ is called the {\em triangulation} of $\Gamma$.
The degree of $\Gamma=|G|$ will be defined as the degree of the combinatorial graph $G$ (the so-defined degree of $\Gamma$ does not depend on the choice of a triangulation). 

It turns out that any graph by a small deformation can be transformed into a graph of degree $\le 3$.

\begin{lemma}\label{degree} For any open cover $\U$ of a topological graph $\Gamma$ there is a $\U$-map $f:\Gamma\to G$ onto a topological graph $G$ of degree $\le 3$.
\end{lemma}

This lemma can be easily proved by induction (and we suspect that it is known as a folklore). The following drawing illustrates how to decrease the degree of a selected vertex of a graph.

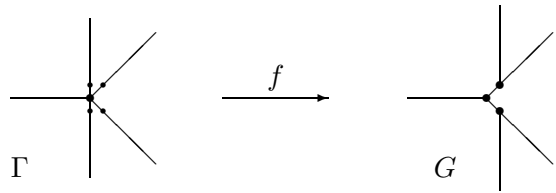
\begin{picture}(100,80)(-120,8)

\put(50,50){\circle*{3}}
\put(50,50){\line(0,1){30}}
\put(50,50){\line(0,-1){30}}
\put(50,50){\line(-1,0){30}}
\put(50,50){\line(1,1){25}}
\put(50,50){\line(1,-1){25}}
\put(55,55){\circle*{2}}
\put(55,45){\circle*{2}}
\put(50,55){\circle*{2}}
\put(50,45){\circle*{2}}
\put(20,20){$\Gamma$}
\put(180,20){$G$}

\put(100,50){\vector(1,0){40}}
\put(117,55){$f$}

\put(200,50){\circle*{3}}
\put(205,55){\line(0,1){30}}
\put(205,45){\line(0,-1){30}}
\put(200,50){\line(-1,0){30}}
\put(200,50){\line(1,1){25}}
\put(200,50){\line(1,-1){25}}
\put(205,55){\circle*{3}}
\put(205,45){\circle*{3}}
\end{picture}

Now we have all tools for the proof of Theorem~\ref{hat}. 
So, take a class $\mathsf K$ of 1-dimensional continua and assume that $X$ is a 4-$\mathsf K$-like continuum. We should prove that $X$ is $\hT$-like.

First, we show that $X$ is 1-dimensional. This will follow from Hemmingsen's Theorem~\ref{hemming} as soon as we check that each open cover $\U$ of $X$ of cardinality $|\U|\le 3$ has an open refinement $\V$ of order $\le 2$. 
Since $|\U|\le 4$ and $X$ is $4$-$\mathsf K$-like, there is a $\U$-map $f:X\to K$ onto a continuum $K\in\mathsf K$. It follows that for some open cover $\V$ of $K$ the cover $f^{-1}(\V)$ refines the cover $\U$. Since the space $K$ is 1-dimensional, the cover $\V$ has an open refinement $\W$ of order $\le 2$. Then the cover $f^{-1}(\W)$ is an open refinement of $\U$ having order $\le 2$.

To prove that $X$ is $\hT$-like, fix any open cover $\U$ of $X$. Because of the compactness of $X$, we can additionally assume that the cover $\U$ is finite. Being 1-dimensional, the continuum $X$ admits a $\U$-map $f:X\to \Gamma$ onto a topological graph $\Gamma$. By Lemma~\ref{degree}, we can assume that $\deg(\Gamma)\le 3$. Adding vertices on edges of $\Gamma$, we can find a triangulation $(V_\Gamma,E_\Gamma)$ of $\Gamma$ so fine that
\begin{itemize}
\item  the path-distance between any vertices of degree 3 in the graph $\Gamma$ is $\ge 6$;
\item the cover $\{f^{-1}(B_2(v))\colon v\in V_\Gamma\}$ of $X$ is inscribed into $\U$.
\end{itemize}

Lemma~\ref{coloring} yields a 4-coloring $\chi:V_\Gamma\to4$ of $V_\Gamma$ such that any two distinct vertices $u,v\in V_\Gamma$ with $d(u,v)\le 2$ have distinct colors. For each color $i\in 4$ consider the open 1-neighborhood $U_i=\bigcup_{v\in\chi^{-1}(i)}B(v)$ of the monochrome set $\chi^{-1}(i)\subseteq V_\Gamma$ in $\Gamma$. Since open 1-balls centered at vertices $v\in V_\Gamma$ cover the graph $\Gamma$, the 4-element family $\{U_i:i\in 4\}$ is an open cover of $\Gamma$.  Then for the 4-element cover $\U_4=\{f^{-1}(U_i)\colon i\in 4\}$ of the $4$-$\mathsf K$-like continuum $X$ we can find a $\U_4$-map $g\colon X\to Y$ to a continuum $Y\in\mathsf K$. Let $\W$ be a finite open cover of $Y$ such that the cover $g^{-1}(\W)$ refines the cover $\U_4$. Since $Y$ is 1-dimensional, we can assume that $\ord(\W)\le 2$. For every $W\in\W$ find a number $\xi(W)\in 4$ such that $g^{-1}(W)\subseteq f^{-1}(U_{\xi(W)})$.

Since $Y$ is a continuum, in particular, a normal Hausdorff space, we may find  a partition of unity
subordinated to the cover $\W$. This is a family $\{\lambda_W\colon W\in\W\}$ of continuous functions $\lambda_W\colon Y\to [0,1]$ such that
\begin{enumerate}
\item[(a)] $\lambda_W(y) = 0$ for $y\in Y \setminus W$,
\item[(b)] $\sum_{W\in\W}\lambda_W(y) = 1$ for all $y\in Y$.
\end{enumerate}

For every $W\in\W$ consider the ``vertical'' family of rectangles 
$$\mathcal R_W=\{W\times B(v):v\in V_\Gamma,\;\chi(v)=\xi(W)\}$$ in $Y\times\Gamma$ and let $\mathcal R=\bigcup_{W\in\W}\mathcal R_W$. For every rectangle $R\in\mathcal R$ choose a set $W_R\in\W$ and a vertex $v_R\in V_\Gamma$ such that $R=W_R\times B(v_R)$. Also let $\mathcal R_R=\{S\in\mathcal R:R\cap S\ne\emptyset\}$. 

\begin{claim}\label{cl1} For any rectangle $R\in\mathcal R$ and a point $y\in W_R$ the set $\mathcal R_{R,y}=\{S\in\mathcal R_R:y\in W_S\}$ contains at most two distinct rectangles.
\end{claim}

\begin{proof} Assume that besides the rectangle $R$ the set $\mathcal R_{R,y}$ contains two other distinct rectangles $S_1=W_{S_1}\times B(v_{S_1})$ and $S_2=W_{S_2}\times B(v_{S_2})$. Taking into account that $y\in W_R\cap W_{S_1}\cap W_{S_2}$ and $\ord(\W)\le 2$, we conclude that either  $W_{S_1}=W_{S_2}$ or $W_R=W_{S_1}$ or $W_R=W_{S_2}$. If $W_{S_1}=W_{S_2}$, then
$$\chi(v_{S_1})=\xi(W_{S_1})=\xi(W_{S_2})=\chi(v_{S_2}).$$
Since $B(v_R)\cap B(v_{S_1})\ne\emptyset\ne B(v_R)\cap B(v_{S_2})$ the property of $4$-coloring $\chi$  implies that $v_{S_1}=v_{S_2}$ and hence $S_1=S_2$. By analogy we can prove that $W_R=W_{S_1}$ implies $R=S_1$ and  $W_R=W_{S_2}$ implies $R=S_2$, which contradicts the choice of $S_1,S_2\in\mathcal R_{R,y}\setminus\{R\}$.   
\end{proof}
 
Claim~\ref{cl1} implies that for every rectangle $R=W_R\times B(v_R)$ the function $\lambda_R:W_R\to \bar B(v_R)\subseteq \Gamma$ defined by
$$\lambda_R(y)=
\begin{cases}\lambda_{W_R}(y)v_R+\lambda_{W_{S}}(y)v_{S},&\mbox{if $\mathcal R_{R,y}=\{R,S\}$ for some $S\ne R$},\\
v_R,&\text{if $\mathcal R_{R,y}=\{R\}$}
\end{cases}
$$
is well-defined and continuous.
Let $\pi_R:R\to W_R\times \bar B(v_R)\subset \bar R$ be the map defined by $\pi_R(y,t)=(y,\lambda_R(y))$.

The graphs of two functions $\lambda_R$ and $\lambda_S$ for two intersecting rectangles $R,S\in\mathcal R$ are drawn on the following picture: 

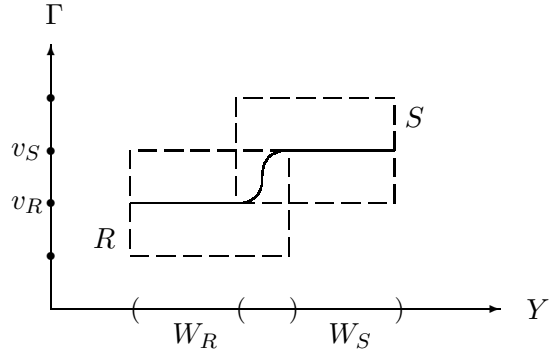
\begin{picture}(200,148)(-120,0)
\put(20,20){\vector(1,0){170}}
\put(200,16){$Y$}
\put(50,17){\small (}
\put(110,17){\small )}
\put(66,7){$W_R$}
\put(90,17){\small (}
\put(150,17){\small )}
\put(125,7){$W_S$}

\put(20,20){\vector(0,1){100}}
\put(18,128){$\Gamma$}

\put(20,40){\circle*{3}}
\put(20,60){\circle*{3}}
\put(5,58){$v_R$}
\put(20,80){\circle*{3}}
\put(5,78){$v_S$}
\put(20,100){\circle*{3}}

\multiput(50,40)(10,0){6}{\line(1,0){7}}
\multiput(110,40)(0,10){4}{\line(0,1){7}}
\multiput(110,80)(-10,0){6}{\line(-1,0){7}}
\multiput(50,80)(0,-10){4}{\line(0,-1){7}}
\put(36,43){$R$}

\multiput(90,60)(10,0){6}{\line(1,0){7}}
\multiput(150,60)(0,10){4}{\line(0,1){7}}
\multiput(150,100)(-10,0){6}{\line(-1,0){7}}
\multiput(90,100)(0,-10){4}{\line(0,-1){7}}
\put(154,89){$S$}

\put(50,60){\line(1,0){40}}
\put(110,80){\line(1,0){40}}
\qbezier(90,60)(100,60)(100,70)
\qbezier(100,70)(100,80)(110,80)
\end{picture}

It follows that for any rectangles $R,S\in\mathcal R$ we get $\pi_R{\restriction} R\cap S=\pi_S{\restriction} R\cap S$, which implies that the union $\pi=\bigcup_{R\in\mathcal R}\pi_R\colon\bigcup\mathcal R\to \bigcup\mathcal R$ is a well-defined continuous function. It is easy to check that for every rectangle $R=W\times B(v)\in\mathcal R$ we get
$$\pi^{-1}(W\times B(v))\subseteq W\times B_2(v).$$

Consider the diagonal product $g\triangle f\colon X\to Y\times\Gamma$. It is easy to check that $(g\triangle f)(X)\subseteq\bigcup\mathcal R$, which implies that the composition $h=\pi\circ (g\triangle f)\colon X\to \bigcup\mathcal R$ is well-defined. We claim that $h$ is a $\U$-map onto the continuum $L=h(X)$, which belongs to the class $\hT$. 

Given any rectangle $R=W\times B(v)\in\mathcal R$, observe that
\begin{multline*}
h^{-1}(R)=(g\triangle f)^{-1}(\pi^{-1}(W\times B(v)))\subseteq (g\triangle f)^{-1}(W\times B_2(v))=\\=g^{-1}(W)\cap f^{-1}(B_2(v))\subseteq f^{-1}(B_2(v))\subseteq U,
\end{multline*}
for some $U\in\U$. Hence $h$ is a $\U$-map.

The projection $\pr_Y:L\to Y$ is locally injective because $L\subseteq \bigcup\mathcal R$ and for every $R\in\mathcal R$ the restriction $\pr_Y{\restriction} R\cap L:R\cap L\to Y$ is injective. Taking into account that $Y\in\mathsf K$, we conclude that  $L\in\hT$, by the definition of the class $\hT$.

\section{Locally injective maps onto tree-like continua and circle}

The following theorem is known for metrizable continua~\cite{Heath}. 

\begin{proposition}\label{li} Each locally injective map $f:X\to Y$ from a continuum $X$ onto a tree-like continuum $Y$ is a homeomorphism.
\end{proposition}

\begin{proof} By the local injectivity of $f$, there is an open cover $\U'$ such that for every $U\in\U'$ the restriction $f{\restriction} U$ is injective.
Let $\U$ be an open cover of $X$ whose second star $\St^2(\U)$ refines the cover $\U'$. Here $\St(U,\U)=\bigcup\{U'\in\U:U\cap U'\ne\emptyset\}$, $\St(\U)=\{\St(U,\U):U\in\U\}$ and $\St^2(\U)=\{\St(U,\St(\U)):U\in\U\}$.

For every $x\in X$ choose a set $U_x\in\U$ that contains $x$. Observe that for distinct points $x,x'\in X$ with $f(x)=f(x')$ the sets $U_x,U_{x'}$ are disjoint. In the opposite case $x,x'\in U_x\cup U_{x'}\subseteq\St(U_x,\U)\subseteq U$ for some set $U\in\U'$, which is not possible as $f{\restriction} U$ is injective. 

Hence for every $y\in Y$ the family $\U_y=\{U_x:x\in f^{-1}(y)\}$ is disjoint.
Since $f$ is closed and surjective, the set $V_y=Y\setminus f(X\setminus \bigcup\U_y)$ is an open neighborhood of $y$ in $Y$ such that $f^{-1}(V_y)\subseteq\bigcup\U_y$. 

Since the continuum $Y$ is a tree-like, the cover $\V=\{V_y:y\in Y\}$ has a finite tree-like refinement $\W$. 
For every $W\in\W$ find a point $y_W\in Y$ with $W\subseteq V_{y_W}$ and consider the disjoint family $\U_W=\{U\cap f^{-1}(W):U\in\U_{y_W}\}$. It follows that $f^{-1}(W)=\bigcup\U_W$ and hence
 $\U_\W=\bigcup_{W\in\W}\U_W$ is an open cover of $X$.

Now we are able to show that the map $f$ is injective. Assuming the converse, find a point $y\in Y$ and two distinct points $a,b\in f^{-1}(y)$. Since $X$ is connected, there is a chain of sets $\{G_1,G_2,\dots,G_n\}\subseteq \U_\W$ such that $a\in G_1$ and $b\in G_n$.
We can assume that the length $n$ of this chain is the smallest possible.
In this case all sets $G_1,\dots,G_n$ are pairwise distinct.

Let us show that $n\ge 3$. In the opposite case $a\in G_1=U_1\cap f^{-1}(W_1)\in\U_\W$, $b\in G_2= U_2\cap f^{-1}(W_2)\in\U_\W$
and $G_1\cap G_2\ne\emptyset$. So, $a,b\in U_1\cup U_2\subseteq\St(U_1,\U)\subseteq U$ for some $U\in\U'$ and then the restriction $f\restriction U$ is not injective.
So $n\ge 3$.

For every $i\le n$ consider the point $y_i=y_{W_i}$ and find sets $W_i\in\W$ and $U_i\in\U_{y_i}$ such that $G_i=U_i\cap f^{-1}(W_i)\in\U_{W_i}$. Then $(W_1,\dots,W_n)$ is a sequence of elements of the tree-like cover $\W$ such that $y\in W_1\cap W_n$ and $W_i\cap W_{i+1}\ne\emptyset$ for all $i<n$. Since the tree-like cover $\W$ does not contain circle-like subfamilies of length $\ge 3$ there are two numbers $1\le i<j\le n$ such that $W_i\cap W_j\ne \emptyset$, $|j-i|>1$ and $\{i,j\}\ne \{1,n\}$. We can assume that the difference $k=j-i$ is the smallest possible. In this case $k=2$. Otherwise, $W_i,W_{i+1},\dots,W_j$ is a circle-like subfamily of length $\ge 3$ in $\W$, which is forbidden. Therefore, $j=i+2$ and the family $\{W_i,W_{i+1},W_{i+2}\}$ contains at most two distinct sets (in the opposite this family is circle-like, which is forbidden). If $W_i=W_{i+1}$, then $U_i=U_{i+1}$ as the family $\U_{W_i}$ is disjoint. The assumption $W_{i+1}=W_{i+2}$ leads to a similar contradiction. It remains to consider the case $W_i=W_{i+2}\ne W_{i+1}$. Since the sets $U_i,U_{i+2}\in\U_{y_i}$ are distinct, there are distinct  points $x_i,x_{i+2}\in f^{-1}(y_i)$ such that $x_i\in U_i$ and $x_{i+2}\in U_{i+2}$. Since $x_i,x_{i+2}\in U_i\cup U_{i+2}\subset\St^2(U_i,\U)\subset U$ for some $U\in\U'$, the restriction $f|U$ is not injective. This contradiction completes the proof.
\end{proof}

\begin{proposition}\label{lii} If  $f:X\to S^1$ is a locally injective map from a continuum $X$ onto the circle $S^1$,  then $X$ is an arc or a circle.
\end{proposition}

\begin{proof} The compact space $X$ has a finite cover by compact subsets that embed into the circle. Consequently, $X$ is metrizable and 1-dimensional.  We claim that $X$ is locally connected. Assuming the converse and applying Theorem 1 of \cite[\S49.VI]{Kur} (or \cite[5.22(b) and 5.12]{Nad}), we could find a convergence continuum $K\subseteq X$. This a non-trivial continuum $K$, which is the limit of a sequence of continua  $(K_n)_{n\in\w}$ that lie in $X\setminus K$.

By the local injectivity of $f$, the continuum $K$ meets some open set $U\subseteq X$ such that $f\restriction U:U\to S^1$ is a topological embedding. The intersection $U\cap K$, being a non-empty open subset of the continuum $K$ is not zero-dimensional. Consequently, its image $f(U\cap K)\subseteq S^1$ also is not zero-dimensional and hence contains a non-empty open subset $V$ of $S^1$. Choose any  point $x\in U\cap K$ with $f(x)\in V$. The convergence $K_n\to K$, implies the existence of a sequence of points $x_n\in K_n$,  $n\in\w$, that converge to $x$. By the continuity of $f$, the sequence $(f(x_n))_{n\in\w}$ converges to $f(x)\in V$. So, we can find a number $n$  such that $f(x_n)\in V\subseteq f(U\cap K)$ and $x_n\in U$. The injectivity of $f\restriction U$ guarantees that $x_n\in U\cap K$ which is not possible as $x_n\in K_n\subset X\setminus K$.

Therefore, the continuum $X$ is locally connected. By the local injectivity, each point $x\in X$ has an open connected neighborhood $V$ homeomorphic to a (connected) subset of $S^1$. Now we see that the space $X$ is a compact 1-dimensional manifold (possibly with boundary). So, $X$ is homeomorphic either to the arc or to the circle. 
\end{proof}

\section{Proof of Theorem~\ref{metr}}

In the proof we shall use the technique of inverse spectra described in \cite[\S2.5]{En} or \cite[Ch.1]{Chi}. Given a continuum $X$ embed it into a Tychonov cube $[0,1]^\kappa$ of weight $\kappa\ge\aleph_0$. 

 Let $A$ be the set of all countable subsets of $\kappa$, partially ordered by the inclusion relation: $\alpha\le\beta$ iff $\alpha\subseteq\beta$.
For a countable subset $\alpha\subseteq\kappa$ let $X_\alpha=\pr_\alpha(X)$ be the projection of $X$ onto the face $[0,1]^\alpha$ of the cube $[0,1]^\kappa$ and $p_\alpha:X\to X_\alpha$ be the projection map. For any countable subsets $\alpha\subseteq\beta$ of $\kappa$ let $p^\beta_\alpha:X_\beta\to X_\alpha$ be the restriction of the natural projection $[0,1]^\beta\to[0,1]^\alpha$. In such a way we have defined an inverse spectrum $\mathcal S=\{X_\alpha,p_\alpha^\beta:\alpha,\beta\in A\}$ over the index set $A$, which is $\w$-complete in the sense that any countable subset $B\subset A$ has the smallest upper bound $\sup B=\bigcup B$ and for any increasing sequence  $\{\alpha_i\}_{i\in\w}\subset A$ with supremum $\alpha=\bigcup_{i\in\w}\alpha_i$ the space $X_\alpha$ is the limit of the inverse sequence $\{X_{\alpha_i}, p_{\alpha_i}^{\alpha_{i+1}},\omega\}$. The spectrum $\mathcal S$ consists of metrizable compacta $X_\alpha$, $\alpha\in A$, and its inverse limit $\varprojlim \mathcal S$ can be identified with the  space $X$.  By Corollary 1.3.2 of \cite{Chi}, the spectrum $\mathcal S$ is factorizing in the sense that any continuous map $f:X\to Y$ to a second countable space $Y$ can be written as the composition $f=f_\alpha\circ p_\alpha$ for some index $\alpha\in A$ and some continuous map $f_\alpha:X_\alpha\to Y$.

Now we are able to prove the ``if'' and ``only if'' parts of Theorem~\ref{metr}.
To prove the ``if'' part, assume that each map $f:X\to Y$ factorizes through a metrizable $n$-$\mathsf K$-like continuum. To show that $X$ is $n$-$\mathsf K$-like, fix any open cover $\U=\{U_1,\dots,U_n\}$ of $X$. By Lemma 5.1.6 of \cite{En}, there is a closed cover $\{F_1,\dots,F_n\}$ of $X$ such that $F_i\subset U_i$ for all $i\le n$. Since $F_i$ and $X\setminus U_i$ are disjoint closed subsets of the compact space $X=\varprojlim \mathcal S$, there is an index $\alpha\in A$ such that for every $i\le n$ the images $p_\alpha(X\setminus U_i)$ and $p_\alpha(F_i)$ are disjoint and hence $W_i=X_\alpha\setminus p_\alpha(X\setminus U_\alpha)$ is an open neighborhood of $p_\alpha(F_i)$. Then $\{W_1,\dots,W_n\}$ is an open cover of $X_\alpha$ such that $p_\alpha^{-1}(W_i)\subseteq U_i$ for all $i\le n$. 

 By our assumption the projection $p_\alpha:X\to X_\alpha$ can be written as the composition $p_\alpha=g\circ\pi$ of a map $\pi:X\to Z$ onto a metrizable $n$-$\mathsf K$-like continuum $Z$ and a map $g:Z\to X_\alpha$.
For every $i\le n$ consider the open subset $V_i=g^{-1}(W_i)$ of $Z$. Since $Z$ is $n$-$\mathsf K$-like, for the open cover $\V=\{V_1,\dots,V_n\}$ of $Z$ there is a $\V$-map $h:Z\to K$ onto a space $K\in\mathsf K$. Then the composition $h\circ\pi:X\to K$ is a $\U$-map of $X$ onto the space $K\in\mathsf K$ witnessing that $X$ is an $n$-$\mathsf K$-like continuum.

Now we shall prove the ``only if'' part of the theorem. Assume that the continuum $X$ is $n$-$\mathsf K$-like.
We shall need the following lemma.

\begin{lemma}\label{l:factor} For any index $\alpha\in A$ there is an index $\beta\ge \alpha$ in $A$ such that for any open cover $\V=\{V_1,\dots,V_n\}$ of $X_\alpha$ there is a map $f:X_\beta\to K$ onto a space $K\in\mathsf K$ such that $f\circ p_\beta:X\to K$ is a $p_\alpha^{-1}(\V)$-map.
\end{lemma}

\begin{proof} Let $\mathcal B$ be a countable base of the topology of the compact metrizable space $X_\alpha$ such that $\mathcal B$ is closed under unions. Denote by $\mathfrak U$ the family of all possible $n$-set covers
$\{B_1,\dots,B_n\}\subseteq\mathcal B$ of $X_\alpha$. It is clear that the family $\mathfrak U$ is countable.

Each cover $\U=\{B_1,\dots,B_n\}\in\mathfrak U$ induces the open cover $p_\alpha^{-1}(\U)=\{p_\alpha^{-1}(B_i):1\le i\le n\}$ of $X$. Since the continuum $X$ is $n$-$\mathsf K$-like, there is a $p_\alpha^{-1}(\U)$-map $f_\U:X\to K_\U$ onto a space $K_\U\in\mathsf K$. By the metrizability of $K_\U$ and the factorizing property of the spectrum $\mathcal S$, for some index $\alpha_\U\ge \alpha$ in $A$ there is  a map $f_{\alpha_U}:X_{\alpha_\U}\to K_\U$ such that $f_\U=f_{\alpha_\U}\circ p_{\alpha_\U}$. Consider the countable set $\beta=\bigcup_{\U\in\mathfrak U}\alpha_\U$, which is the smallest lower bound of the set $\{\alpha_\U:\U\in\mathfrak U\}$ in $A$. We claim that this index $\beta$ has the required property.

Let $\V=\{V_1,\dots,V_n\}$ be any open cover of $X_\alpha$. By Lemma~5.1.6 of \cite{En}, there is a closed cover $\{F_1,\dots,F_n\}$ of $X_\alpha$ such that $F_i\subseteq V_i$ for all $i\le n$. Since $\mathcal B$ is the base of the topology of $X_\alpha$ and $\mathcal B$ is closed under finite unions, for every $i\le n$ there is a basic set $B_i\in\mathcal B$ such that $F_i\subseteq B_i\subseteq V_i$. Then the cover $\U=\{B_1,\dots,B_n\}$ belongs to the family $\mathfrak U$ and refines the cover $\V$. Consider the map  $f=f_{\alpha_\U}\circ p^\beta_{\alpha_\U}:X_\beta\to K$ and observe that $f\circ p_\beta=f_{\alpha_\U}\circ p_{\alpha_\U}$ is a $p_\alpha^{-1}(\U)$-map and a $p_\alpha^{-1}(\V)$-map.
\end{proof}

Now let us return back to the proof of the theorem. Given a map $f:X\to Y$ to a second countable space, we need to find a map $\pi:X\to Z$ onto a metrizable $n$-$\mathsf K$-like continuum $Z$ and a map $g:Z\to Y$ such that $f=g\circ\pi$. Since the spectrum $\mathcal S$ is factorizing, there are an index $\alpha_0\in A$ and a map $f_0:X_{\alpha_0}\to Y$ such that $f=f_0\circ p_{\alpha_0}$. Using Lemma~\ref{l:factor}, by induction construct an increasing sequence $(\alpha_n)_{n\in\w}$ in $A$ such that for every $i\in\w$ and any open cover $\V=\{V_1,\dots,V_n\}$ of  $X_{\alpha_i}$ there is a map $f:X_{\alpha_{i+1}}\to K$ onto a space $K\in\mathsf K$ such that $f\circ p_{\alpha_{i+1}}$ is a $p_{\alpha_{i}}^{-1}(\V)$-map. 

Let $\alpha=\sup_{i\in\w}\alpha_i=\bigcup_{i\in\w}\alpha_i$. We claim that the metrizable continuum $X_\alpha$ is $n$-$\mathsf K$-like. Given any open cover $\U=\{U_1,\dots,U_n\}$ of $X_\alpha=\varprojlim X_{\alpha_i}$, we can find $i\in\w$ such that the sets $W_i=X_{\alpha_i}\setminus p^\alpha_{\alpha_i}(X_\alpha\setminus U_i)$, $i\le n$, form an open cover $\W=\{W_1,\dots,W_n\}$ of $X_{\alpha_i}$ such that the cover $(p^\alpha_{\alpha_i})^{-1}(\W)$ refines the cover $\U$. By the choice of the index $\alpha_{i+1}$, there is a map $g:X_{\alpha_{i+1}}\to K$ onto a space $K\in\mathsf K$ such that $g\circ p_{\alpha_{i+1}}:X\to K$ is a $p_{\alpha_{i}}^{-1}(\W)$-map. It follows that $g\circ p_{\alpha_{i+1}}^\alpha:X_\alpha\to K$ is a $(p^\alpha_{\alpha_i})^{-1}(\W)$-map and hence a $\U$-map, witnessing that the continuum $X_\alpha$ is $n$-$\mathsf K$-like.

Now we see that the metrizable $n$-$\mathsf K$-like continuum $X_\alpha$ and the maps $\pi=p_\alpha:X\to X_\alpha$ and $g=f_0\circ p^\alpha_{\alpha_0}:X_\alpha\to Y$ satisfy our requirements.

\section{Proof of Theorem~\ref{main}}

The ``chainable and tree-like'' parts of Theorem~\ref{main} follow immediately from the characterization Theorem~\ref{Th3}. So, it remains to prove the ``circle-like'' part. Let $\mathsf K=\{S^1\}$. We need to prove that each 4-$\mathsf K$-like continuum $X$ is $\mathsf K$-like. Given an open cover $\U$ of $X$ we need to construct a $\U$-map of $X$ onto the circle. By Theorem~\ref{metr}, there is a $\U$-map onto a metrizable 4-$\mathsf K$-like continuum $Y$. It follows that for some open cover $\V$ of $Y$ the cover  $f^{-1}(\V)$ refines $\U$. The proof will be complete as soon as we prove that the continuum $Y$ is circle-like. In this case there is a $\V$-map $g:Y\to S^1$ and the composition $g\circ f:X\to S^1$ is a required $\U$-map witnessing that $X$ is circle-like.

By Theorem~\ref{hat}, the metrizable continuum $Y$ is $\hT$-like. By Proposition~\ref{lii}, each continuum $K\in\hT$ is homeomorphic to $S^1$ or $[0,1]$. Consequently, the continuum $Y$ is circle-like or chainable. In the first case we are done. So, we assume that $Y$ is chainable. 

By  \cite[Theorem 12.5]{Nad}, the continuum $Y$ is irreducible between some points $p,q\in Y$.
The latter means that each subcontinuum of $X$ that contains the points $p,q$ coincides with $Y$.
 We claim that $Y$ is either indecomposable or $Y$ is the union of two indecomposable subcontinua. For the proof of this fact we will use the argument of~\cite[Exercise~12.50]{Nad} (cf. also~\cite[Theorem~3.3]{Kra}). 

Suppose that $Y$ is not indecomposable. It means that there are two proper subcontinua $A,B$ of $Y$ such that $Y=A\cup B$.  By the choice of the points $p,q$, they cannot simultaneously lie in $A$ or in $B$. So, we can assume that $p\in A$ and $q\in B$. 

We claim that the closure of the set $Y\setminus A$ is connected. Assuming that $\overline{Y\setminus A}$ is disconnected, we can find a proper closed-and open subset $F\subsetneq\overline{Y\setminus A}$ that contains the point $q$ and conclude that $F\cup A$ is a proper subcontinuum of $Y$ that contains both points $p,q$, which is not possible. Replacing $B$ by the closure of $Y\setminus A$, we can assume that $Y\setminus A$ is dense in $B$.
Then $Y\setminus B$ is dense in $A$.

We claim that the sets $A$ and $B$ are indecomposable. Assuming that $A$ is decomposable, find two proper subcontinua $C,D$ such that $C\cup D=A$. 
We can assume that $p\in D$. Then $B\cap D=\emptyset$ (as $Y$ is irreducible between $p$ and $q$). By Theorem~11.8 of~\cite{Nad}, the set $Y\setminus(B\cup D)$ is connected. Let $Z_1$ and $Z_2$ be open disjoint subsets of $X$ such that $B\subseteq Z_1$ and $D\subseteq Z_2$. Since $Y$ is $4$-$\{S^1\}$-like, for the open cover $\mathcal{Z}=\{Z_1,Z_2,Y\setminus(B\cup D)\}$ of $Y$ there exists a $\mathcal{Z}$-map $h\colon Y\to S^1$. Thus $h(B)\cap h(D)=\emptyset$ and $S^1\setminus (h(B)\cup h(D))$ is the union of two disjoint open intervals $W_1,W_2$. Since $h$ is a $\mathcal{Z}$-map, $Y\setminus(B\cup D)= h^{-1}(W_1)\cup h^{-1}(W_2)$ which contradicts the connectedness of the set $Y\setminus (B\cup D)$.

Now we know that $Y$ is either indecomposable or is the union of two indecomposable subcontinua. Applying Theorem 7 of~\cite{Bur}, we conclude that the metrizable chainable continuum $Y$ is circle-like. 

\section{Open Problems}

\begin{problem} For which families $\mathsf{K}$ of connected topological graphs every $4$-$\mathsf{K}$-like continuum is $\mathsf{K}$-like? Is it true for the family $\mathsf K=\{8\}$ that contains 8, the bouquet of two circles?
\end{problem}

Also we do not know if Theorem~\ref{Th3} can be generalized to classes of higher-dimensional continua.

\begin{problem} Let $k\in\IN$ and $\mathsf K$ be a class of $k$-dimensional (contractible) continua. Is there a finite number $n$ such that a continuum $X$ is $\mathsf K$-like if and only if it is $n$-$\mathsf K$-like?
\end{problem}

\section*{Acknowledgment} 
The authors express their sincere thanks to the anonymous referee whose valuable remarks and suggestions helped the authors to improve substantially the results of this paper.

\end{document}